\magnification \magstep 1 
\centerline{\bf Presentations for subgroups of Artin groups} 
\medskip
\centerline{by Warren Dicks and Ian J. Leary} 
\medskip
\def\bb{1}
\def\bie{2} 
\def\chis{3}
\def\dic{4}
\def\dro{5}
\def\jho{6}
\def\kir{7} 
\def\sta{8}
\def\Bbb#1{{\bf #1}}
\def\ignore#1{}
\def\From{From}
\def\pra{\par}
\def\proof{\par \noindent Proof.\quad}
\def\hhrule#1#2{\kern-#1
   \hrule height#1 depth#2 \kern-#2 }
\def\hvrule#1#2{\kern-#1{\dimen0=#1
    \advance\dimen0 by#2\vrule width\dimen0}\kern-#2 }
\def\makeblankbox#1#2{\setbox0=\hbox{A}
\hbox{\lower\dp0\vbox{\hhrule{#1}{#2}%
   \kern -#1
   \hbox to \wd0{\hvrule{#1}{#2}%
     \raise\ht0\vbox to #1{}
     \lower\dp0\vtop to #1{}
     \hfil\hvrule{#2}{#1}}
  \kern-#1\hhrule{#2}{#1}}}}
\def\qed{\makeblankbox{0pt}{.3pt}}
\def\remark{\par\medskip\noindent {\bf Remark.}\quad} 

\noindent 
Let $\Delta$ be a finite flag complex, that is, a finite
 simplicial complex that
contains a simplex bounding every complete subgraph of its
1-skeleton.  The associated right-angled Artin group $G_\Delta$ is the
group given by the presentation with generators the vertex set of
$\Delta$, and relators the commutators $[v,w]$ for each pair of
adjacent vertices in $\Delta$.  For example, the $n$-simplex
corresponds to a free abelian group of rank $n+1$, a complex
consisting of $n$ points corresponds to the free group $F_n$ of rank
$n$, the group corresponding to the square is $(F_2)^2$, and the group
corresponding to the octahedron is $(F_2)^3$.  

Provided that $\Delta$ is non-empty, there is a homomorphism from
$G_\Delta$ onto the integers, that takes every generator to 1.  The
group $H_\Delta$ is defined to be the kernel of this homomorphism.  
Remarkable recent work of Mladen Bestvina and Noel Brady has shown
that the homological finiteness properties of the group $H_\Delta$ are 
controlled by the topology of the complex $\Delta$ [\bb].  They show that 
$H_\Delta$ is finitely generated if and only if $\Delta$ is connected,
$H_\Delta$ is finitely presented if and only if $\Delta$ is
1-connected, and $H_\Delta$ is of type $FP(n)$ if and only if $\Delta$
is $(n-1)$-acyclic.  Precursors of this result include J. Stallings'
group that is finitely presented but not of type $FP(3)$, which is 
$H_\Delta$ in the case when $\Delta$ is the octahedron [\sta], 
and R. Bieri's group of type $FP(n)$ but not of type $FP(n+1)$, which
is $H_\Delta$ in the case when $\Delta$ is a join of $(n+1)$ pairs of 
points [\bie,2.6].  

The arguments used by Bestvina and Brady are geometric, and they do
not give presentations for the groups that they consider.  Theorem 1
of this paper 
gives a presentation for $H_\Delta$ for any connected $\Delta$.  The
generators in the presentation are the edges of $\Delta$, and each
1-cycle in $\Delta$ gives rise to an infinite family of relators. 
In the case when $\Delta$ is simply connected, it is shown how to
reduce this presentation to a finite one.  This gives an independent
and purely algebraic proof that $H_\Delta$ is finitely presented
when $\Delta$ is simply connected.  It would be interesting to give a 
similar proof of the converse.  In Proposition 4 and Corollary 5, we 
review some results concerning the homology of $G_\Delta$.  Using
an argument due to Stallings, we deduce that when
the Euler characteristic of $\Delta$ is not equal to that of a point, the
rational cohomology of $H_\Delta$ cannot be finite dimensional.

\proclaim Definition.  For $e$ a directed edge of $\Delta$, let
$\iota e$ (resp.\ $\tau e$) denote the initial (resp.\ terminal)
vertex of $e$.  

\proclaim Theorem 1.  Let $\Delta$ be connected.  The group
$H_\Delta$ has a presentation with generators the set of directed
edges of $\Delta$, and relators all words of the form
$e_1^ne_2^n\ldots e_l^n$, where $l,n\in \Bbb Z$, $n\neq 0$, $l\geq 2$,
and $(e_1,\ldots,e_l)$ is a directed cycle in $\Delta$.  In terms of the
given generators for $G_\Delta$, $e= \iota e(\tau e)^{-1}$.

\proof  Let $H'_\Delta$ be the group presented as in the statement.  
Cycles of length two in $\Delta$ have the form $(e,\bar e)$, where $\bar
e$ is the same edge as $e$ with its opposite orientation.  The
relation $e\bar e=1$ implies that for each edge $e$, $\bar e = e^{-1}$
in $H'_\Delta$.  These relations could of course be used to halve the
size of the generating set at the expense of adding some signs in the 
relators.  We shall apply these relations as necessary without comment
in the sequel.  \From\ the symmetrical form of the relators in
$H'_\Delta$, it follows that there is an endomorphism $\xi$ of
$H'_\Delta$ which sends each directed edge $e$ to $e^{-1}$.  Moreover,
$\xi^2=1$, and so $\xi$ is an automorphism.  

Define a
homomorphism $\phi$ from $H'_\Delta$ to $G_\Delta$ by $\phi(e)=
\iota e(\tau e)^{-1}$.  This does define a homomorphism, since if
$(e_1,\ldots,e_l)$ is a directed cycle, with $\iota e=a_i$ and 
$\tau e=a_{i+1}$, then 
$$\phi(e_1^n\cdots e_l^n) = (a_1^na_2^{-n})(a_2^na_3^{-n})\cdots 
(a_l^na_1^{-n}) = 1.$$
The image of $\phi$ is contained in $H_\Delta$.  Conversely, any
element $w$ of $H_\Delta$ is expressible in the form
$w=a_1^{n(1)}a_2^{n(2)}\cdots a_m^{n(m)}$, where $n(1)+\cdots +n(m)=0$.
If $(e_1,\ldots,e_r)$ is a directed path from $a_m$ to $a_{m-1}$, then 
$\phi(e_1^{-n(m)}\cdots e_r^{-n(m)}) = a_m^{-n(m)}a_{m-1}^{n(m)}$.  
By induction on $m$ it follows that any element of $H_\Delta$ is in
the image of $\phi$.  

It remains to prove that $\phi$ is injective.  To do this, our
strategy is as follows.  Define an extension $G'_\Delta$ with kernel
$H'_\Delta $ and quotient $\Bbb Z$, extend $\phi$ to a homomorphism 
$\tilde\phi$ from $G'_\Delta$ to $G_\Delta$, and finally show that 
$\tilde\phi$ is an isomorphism by exhibiting an inverse.  

For vertices $a$ and $b$ of $\Delta$, define $p(a,b)$ to be element of
$H'_\Delta$ represented by the word $e_1\cdots e_n$, where
$(e_1,\ldots,e_n)$ is a choice of directed edge path from $a$ to $b$.
The hypothesis that $\Delta$ be connected ensures that there is such a
choice.  Note that different choices of directed edge path give rise
to the same element of $H'_\Delta$.  Fix a vertex $a$, and define an
endomorphism $\psi_a$ of $H'_\Delta$ by 
$\psi_a(e)=p(a,\iota e) e p(\iota e,a)$.  It must be shown that this
does define an endomorphism.  First note that 
$$\psi_a(e^n)= p(a,\iota e) e^n p(\iota e,a) = p(a,\iota e)e^{n+1} 
p(\tau e, a).$$
Now if $(e_1,\ldots,e_l)$ is any directed path from $b$ to $b'$, then
for any non-zero integer $n$,
$$\psi_a(e_1^ne_2^n\cdots e_l^n) = p(a,\iota e_1)e_1^{n+1}
p(\tau e_1,a)p(a,\iota e_2)e_2^{n+1}p(\tau e_2,a)\cdots p(a,\iota e_l)
e_l^{n+1}p(\tau e_l,a),$$
which telescopes to the equation  
$$\psi_a(e_1^n\cdots e_l^n)= 
p(a,b)e_1^{n+1}e_2^{n+1}\cdots e_l^{n+1}p(b',a).$$
By considering the case when $b=b'$, it may be seen that $\psi_a$
defines an endomorphism of $H'_\Delta$.  
 
To show that $\psi_a$ is an automorphism, we compute its composite
with the automorphism $\xi$ defined by $\xi(e)= e^{-1}$.  The case
$n=-1$ of the equation above gives 
$$\psi_a\xi(p(b,b'))= p(a,b)p(b',a).$$
Hence for any edge $e$, 
$$\eqalign{\psi_a\xi\psi_a\xi(e)
&= \psi_a\xi\left(p(a,\iota e)p(\tau e,a)\right)\cr 
&= \psi_a\xi\left(p(a,\iota e)\right) \psi_a\xi\left(p(\tau
e,a)\right)\cr 
&= \left(p(a,a)p(\iota e,a)\right)\left(p(a,\tau e)p(a,a)\right) \cr 
&= e.\cr}$$ 
Thus $\psi_a\xi$ has order two, so is an automorphism, and it follows
that $\psi_a$ is an automorphism.

Let $G'_\Delta$ be the extension with kernel $H'_\Delta$ and infinite
cyclic quotient, generated by $a'$, where the conjugation action of
$a'$ is given by $a'ha'^{-1}= \psi_a(h)$.  For every directed edge $e$,
$\phi\left(\psi_a(e)\right) = a\phi(e)a^{-1}$, and so $\phi$ may be
extended to a surjective homomorphism $\tilde \phi \colon G'_\Delta
\rightarrow G_\Delta$ by setting $\tilde \phi(a')= a$.  Define a
map $\theta\colon G_\Delta \rightarrow G'_\Delta$ on the generators 
by $\theta(b) = p(b,a)a'$.  For any directed edge $e$, $p(\iota e, a)
= ep(\tau e,a)$ in $H'_\Delta$, and hence $\theta(\iota e)=
e\theta(\tau e)$.  It is also true that 
$$\eqalign{\theta(\tau e)e &= p(\tau e,a)a'e\cr 
&= p(\tau e,a)\psi_a(e)a' \cr 
&=p(\tau e,a)p(a,\iota e)e p(\iota e,a)a' \cr
&= p(\iota e,a)a' = \theta(\iota e).\cr}$$ 
Hence for each directed edge $e$, $\theta(\tau e)^{-1}\theta(\iota e) = e = 
\theta(\iota e)\theta(\tau e)^{-1}$. 
It follows that $\theta$ extends to a homomorphism from $G_\Delta$ to
$G'_\Delta$.  Finally, it may be shown that $\theta$ is an inverse to 
$\tilde\phi$, which is therefore an isomorphism.  
\qed

\proclaim Proposition 2.  For a directed cycle $c= (e_1,\ldots,e_l)$ 
of length $l$ in $\Delta$,
write $c^{[n]}$ for the relator $e_1^n\cdots e_l^n$ in the above
presentation for $H_\Delta$.  The presentation may be simplified as
follows.  
\pra 
(a) When $l=2$, the relators $c^{[n]}$ are consequences of $c^{[1]}$.  
\pra 
(b) When $l=3$, the relators $c^{[n]}$ are consequences
of the relators $c^{[1]}$ and $c^{[-1]}$.  
\pra 
(c) Let $c_1,\ldots,c_m$ be directed cycles such that the normal
subgroup of $\pi_1(\Delta)$ generated by the conjugacy classes of 
$c_1,\ldots,c_m$ is the whole of $\pi_1(\Delta)$.  Then each relator 
$c^{[n]}$ is a consequence of the (finitely many) relators given in a)
and b), together with the relators $c_i^{[n]}$ for all $n$ and $1\leq
i\leq m$.  

\proof As remarked earlier, the relation $e\bar e=1$ implies that 
$\bar e = e^{-1}$.  The relations $e^n\bar e^n$ for all $n$ are
consequences of this.  Now suppose that $(e,f,g)$ is a directed triangle
in $\Delta$.  The relation $efg=1$ implies that $g=f^{-1}e^{-1}$.
Substituting for $g$, the relator $e^{-1}f^{-1}g^{-1}$ becomes the
commutator $e^{-1}f^{-1}ef$.  It follows that the group generated by
$e$, $f$ and $g$ subject to just these two relators is isomorphic to 
$\Bbb Z \times \Bbb Z$, via an isomorphism sending $e$, $f$, $g$ to 
$(1,0)$, $(0,1)$ and $(-1,-1)$.  For each $n$, the relation 
$e^nf^ng^n=1$ already holds in this group.  

If the cycles $c$ and $c'$ differ only in their choice of starting
point, i.e., $c=(e_1,\ldots, e_l)$ and $c'=(e_{r+1},\ldots, e_l,e_1,\ldots,
e_r)$, then the relator $c'^{[n]}$ is a conjugate of the relator 
$c^{[n]}$.  Similarly, if $c'$ is equal to $c$ but in the opposite
orientation, then $c'^{[-n]}$ is a conjugate of $c^{[n]}$.  If $c$ is a
composite cycle, $c=c'.c''$, then $c^{[n]}$ is equal to 
$c'^{[n]}c''^{[n]}$.  Let $p$ be any directed edge path in $\Delta$, 
and let $\bar p$ be the opposite path.  The
cycle $\bar p.p$ is a composite of length two cycles of the form
$e\bar e$, so the relators $(\bar pp)^{[n]}$ are consequences of the
finitely many relators listed in (a).   

For convenience, choose a basepoint $a$ for $\Delta$.  For each path 
$p$ from $a$ to the initial point of $c_i$, the relator $(pc_i\bar
p)^{[n]}$ is a consequence of $c_i^{[n]}$ and the relators listed in (a).  
By assumption, composites of based cycles of this form and their
inverses suffice to construct an element of every based homotopy class
of cycles in $\Delta$.  To complete the proof, it suffices to show
that whenever $c$ and $c'$ are homotopic cycles, the relator $c'^{[n]}$
is a consequence of $c^{[n]}$ and the relators listed in (a) and (b).  
View a simplicial homotopy from $c$ to $c'$ as a triangulation of
$S^1\times I$ together with a labelling of its directed edges by
generators of 
$H_\Delta$ and the identity element, having the following properties:
the labelling of $S^1\times \{0\}$ spells $c$; the labelling of $S^1\times
\{1\}$ spells $c'$; the labelling of each triangle is of one of
the following forms:  $(e,f,g)$ for some triangle $(e,f,g)$ in $\Delta$,
or $(e,\bar e,1)$, or $(1,1,1)$.  Now move from $S^1\times \{0\}$ to
$S^1\times \{1\}$ \lq one simplex at a time', at each stage either
removing a 1-simplex that is connected to the remaining complex only
at one end, or removing a 2-simplex that has an edge in the boundary
of the remaining complex, together with that edge.  This gives rise
to a sequence $c=d_0,d_1,\ldots,d_k=c'$ of cycles in $\Delta$ such
that $d_{i+1}$ differs from $d_i$ in only one of the following ways:  
\def\hand{\hbox{and}}
$$\eqalign{d_{i+1} &= d_i, \cr 
d_{i+1} = d'.e\bar. e.d''\quad &\hand\quad d_i =
d'.d'',\quad\hbox{for some edge $e$,} \cr 
d_i = d'.e.\bar e.d''\quad &\hand\quad d_{i+1} = d'.d'',
\quad\hbox{for some edge $e$,} \cr 
d_{i+1} = d'.\bar g.\bar f.d''\quad &\hand\quad d_i = d'.e.d'',\quad 
\hbox{for some triangle $(e,f,g)$.}\cr
}$$ 
In each case it may be seen that the relator $d_{i+1}^{[n]}$ is a
consequence of $d_i^{[n]}$ together with the relators of types (a) and
(b).  
\qed

\proclaim Corollary 3.  When $\Delta $ is simply connected, $H_\Delta$
has a finite presentation.  

\proof This is a special case of Proposition 2(c).  \qed

\remark Let $K_\Delta$ be the group generated by the directed edges of
the complex $\Delta$ subject only to the relations of types (a) and
(b) of Proposition 2.  There is a natural map from $K_\Delta$ onto 
$H_\Delta$, which is an isomorphism when $\Delta$ is simply
connected.  To show that $H_\Delta$ cannot be finitely presented when
$\Delta$ is not simply connected, it would suffice to show that 
$\ker(K_\Delta \rightarrow H_\Delta)$ is not finitely generated as a
normal subgroup of $H_\Delta$.  It would be very interesting if one 
could do this algebraically.  The
group $K_\Delta$ arises implicitly in the work of Mladen Bestvina
and Noel Brady, as the fundamental group of a \lq level set' [\bb].
Jim Howie constructs another 2-complex with fundamental group
$K_\Delta$, and uses this to give an alternative proof that $H_\Delta$
is $FP(2)$ when $\Delta$ is 1-acyclic, in [\jho].

\proclaim Definition.  For a set $V$, let $\Lambda^*(V)$ be the exterior
algebra generated (in degree one) by the functions from $V$ to $\Bbb
Z$.  For a simplicial complex $\Delta$, the exterior
face ring $\Lambda^*_\Delta$ is the quotient $\Lambda(V)/I$, where $V$
is the vertex set of $\Delta$, and $I$
is the ideal generated by all monomials of the form $v_0\cdots v_n$ 
such that $(v_0,\ldots,v_n)$ is not an $n$-simplex of $\Delta$.

The exterior face ring $\Lambda^*_\Delta$ is a contravariant functor of
$\Delta$, and is graded, since $I$ is generated by homogeneous
elements.  For each $i>0$, $\Lambda_\Delta^i$ is a free abelian group
of rank equal to the number of $(i-1)$-simplices in $\Delta$.  
The complex $\Delta$ is a flag complex if and only if $I$
is generated by monomials of degree two.  

\proclaim Proposition 4 (Kim-Roush).  The integral cohomology 
ring of $G_\Delta$ is 
naturally isomorphic to the exterior face ring of $\Delta$.  

\proof Since this proposition has already appeared in [\kir], p.~185
and in [\dic], pp.~227--228, we merely sketch the proof.  When
$\Delta$ is a simplex, $G_\Delta$ is free abelian and 
$H^*(G_\Delta)=\Lambda({\rm Hom}(G,\Bbb Z))= \Lambda_\Delta$.  
When $\Delta$ is not a simplex, $\Delta= \Delta_1\cup \Delta_2$ 
for smaller flag complexes $\Delta_1$ and $\Delta_2$, with 
$\Delta_3=\Delta_1\cap\Delta_2$.  This gives rise to a decomposition
of $G$ as a free product with amalgamation $G=G_1*_{G_3}G_2$, where 
$G_i = G_{\Delta_i}$.  The claim for $\Delta$ follows by a 
Mayer-Vietoris argument.  \qed

The following result is due to C. Droms [\dro], see also 
I. Chiswell's generalization [\chis].   

\proclaim Corollary 5 (Droms).  The Euler characteristic of $G_\Delta$
is equal to $1-\chi(\Delta)$.  
\qed

\proclaim Corollary 6.  If the rational cohomology of $H_\Delta$ is
finite-dimensional, then $\chi(\Delta)=1$.  

\proof Let $G=G_\Delta$ and $H=H_\Delta$.  The Mayer-Vietoris 
sequence for $G$ expressed as an HNN-extension with base group $H$ is
a long exact sequence:   
$$\cdots \rightarrow H^{i-1}(H)\rightarrow H^{i-1}(H)
\rightarrow H^i(G) \rightarrow H^i(H)\rightarrow
H^i(H) \rightarrow \cdots.$$
In the case when $H^*(H)=H^*(H,\Bbb Q)$ is finite-dimensional, it
follows that the Euler characteristic of $G$ must be zero, and hence
that $\chi(\Delta)=1$.   
\qed

\proclaim Corollary 7 (Bestvina-Brady).  If $\Delta$ is connected and 
simply connected, and $\chi(\Delta)\neq 1$, then $H_\Delta$ is
finitely presented but not of type $FP$.  
\qed 

\medskip
\noindent 
{\bf Acknowledgements.} The first author acknowledges support from the
DGICYT (Spain) through grant number PB93-0900, and the second author 
acknowledges support from the Nuffield Foundation through grant number
SCI/180/96/127.  

\bigskip 
\noindent 
{\bf References.} 
\medskip
\par\frenchspacing

\item{[\bb]} M. Bestvina and N. Brady, Morse theory and finiteness
properties of groups, to appear in {\it Invent. Math.}.  

\item{[\bie]} R. Bieri, Homological dimension of discrete groups, {\it
Queen Mary College Mathematics Notes}, University of London (1976).  

\item{[\chis]} I. M. Chiswell, The Euler characteristic of graph
products and of Coxeter groups, {\it Discrete groups and Geometry
(Birmingham 1991)}, London Math. Soc. Lecture Notes {\bf 173}, 36--46,
Cambridge Univ. Press, Cambridge (1992).  

\item{[\dic]} W. Dicks, An exact sequence for rings of polynomials in
partly commuting indeterminates, {\it J. Pure Appl. Algebra} {\bf
22} (1981), 215--228.  

\item{[\dro]} C. Droms, Subgroups of graph groups, {\it J. Algebra}
{\bf 110} (1987), 519--522.  

\item{[\jho]} J. Howie, Bestvina-Brady groups and the plus
construction, preprint (1997).  

\item{[\kir]} K. H. Kim and F. W. Roush, Homology of certain algebras
defined by graphs, {\it J. Pure Appl. Algebra} {\bf 17} (1980),
179--186.  

\item{[\sta]} J. Stallings, A finitely presented group whose
3-dimensional integral homology is not finitely generated, 
{\it Amer. J. Math.} {\bf 85} (1963), 541--543.  

\bigskip
\bigskip
\hbox{
\vbox{
\hbox{Warren Dicks} 
\hbox{Departament de Matem\`atiques,} 
\hbox{Universitat Aut\`onoma de Barcelona,} 
\hbox{E 08193 Bellaterra (Barcelona),} 
\hbox{Spain.}
\hbox{\phantom{United Kingdom}}} 
\qquad 
\vbox{
\hbox{Ian J. Leary}
\hbox{Faculty of Mathematical Studies,} 
\hbox{University of Southampton,} 
\hbox{Southampton,}
\hbox{SO17 1BJ,}
\hbox{United Kingdom.} 
}}

\end